\documentclass[11pt]{amsart}

\usepackage{xspace}

\newcommand{\punctual}{degree~$0$\xspace}

\usepackage{amssymb}

\usepackage{array,dcolumn}
\newcolumntype{C}{>{$}c<{$}}
\newcolumntype{L}{>{$}l<{$}}
\newcolumntype{R}{>{$}r<{$}}

\def\to{\mathchoice
{\longrightarrow}
{\rightarrow}
{\rightarrow}
{\rightarrow}}

\def\mapsto{\DOTSB\mapstochar\to}

\newtheorem{theorem}{Theorem}
\newtheorem*{theorem*}{Theorem}

\newtheorem{lemma}[theorem]{Lemma}
\newtheorem{corollary}[theorem]{Corollary}

\theoremstyle{definition}

\newcommand{\Z}{\mathbb{Z}}
\newcommand{\C}{\mathbb{C}}

\newcommand{\Q}{\mathbb{Q}}

\newcommand{\om}{\omega}

\newcommand{\p}{\partial}
\renewcommand{\*}{\cdot}
\newcommand{\half}{\tfrac12}

\newcommand{\bull}{\bullet}

\newcommand{\<}{\langle}
\renewcommand{\>}{\rangle}

\renewcommand{\[}{{[\![}}

\newcommand{\CM}{\mathcal{M}}

\newcommand{\Mbar}{\overline{\mathcal{M}}}

\DeclareMathOperator{\Vir}{\mathfrak{v}}
\newcommand{\CE}{\mathcal{E}}
\newcommand{\CO}{\mathcal{O}}
\newcommand{\EE}{\mathbb{E}}
\newcommand{\Dual}{\vee}
\renewcommand{\t}{\mathbf{t}}
\newcommand{\s}{\mathbf{s}}

\renewcommand{\a}{\mathbf{a}}
\renewcommand{\b}{\mathbf{b}}
\renewcommand{\c}{\mathbf{c}}
\renewcommand{\P}{\mathbb{P}}

\newcommand{\Alpha}{\boldsymbol{\alpha}}
\newcommand{\Beta}{\boldsymbol{\beta}}

\begin{document}

\title{Virasoro constraints and the Chern classes of the Hodge bundle}

\author{E. Getzler}

\address{Northwestern University, Evanston, IL 60208}

\email{getzler@math.nwu.edu}

\author{R. Pandharipande}

\address{University of Chicago, Chicago, IL 60637}

\email{rahul@math.uchicago.edu}

\thanks{Communications with C. Faber were invaluable in our research.
This paper was written while the authors were guests of the Scuola
Normale Superiore in Pisa. Both authors are partially supported by the NSF}

\maketitle

\section{Introduction}

\subsection{Virasoro constraints}
Let $X$ be smooth projective variety of dimension $r$. The descendent
Gromov-Witten invariants $\< \tau_{k_1}(\gamma_{a_1}) \dots
\tau_{k_n}(\gamma_{a_n}) \>_{g,\beta}^X$ of $X$ may be assembled into
generating functions%
\footnote{The peculiar ordering of the variables $t^{b_i}_{\ell_i}$ is in
case there are odd dimensional cohomology classes on $X$.}
\begin{multline}
\< \< \tau_{k_1}(\gamma_{a_1}) \dots \tau_{k_n}(\gamma_{a_n}) \> \>_g^X \\
= \sum_{\beta \in H_2(X,\Z)} q^\beta \sum_{N\ge 0} \frac{1}{N!}
\sum_{\substack{b_1\dots b_N \\ \ell_1 \dots \ell_N}} t_{\ell_N}^{b_N}
\dots t_{\ell_1}^{b_1} \< \tau_{\ell_1}(\gamma_{b_1}) \dots
\tau_{\ell_N}(\gamma_{b_N}) \tau_{k_1}(\gamma_{a_1}) \dots
\tau_{k_n}(\gamma_{a_n}) \>_{g,\beta}^X .
\end{multline}
In particular, we may form the exponential generating function for all
of the descendent Gromov-Witten invariants
\begin{equation} \label{GW}
Z^X = \exp \biggl( \sum_{g\ge 0} \hbar^{g-1} \<\<~\>\>^X_g \biggr) .
\end{equation}

It has been conjectured by Eguchi, Hori, and Xiong \cite{EHX1,EHX2} that
$Z^X$ is annihilated by formal differential operators $L_k$, $k\ge-1$, on
the affine space with coordinate $\{t_k^a\mid k\ge0\}$, whose definition
(which we recall in Section \ref{Virasoro}) depends only on the inner
product space $H^\bull(X,\C)$, its Hodge decomposition, and the
endomorphism of multiplication by the anticanonical class $c_1(X)$, and
which satisfy the commutation relations
$$
[L_k,L_\ell] = (k-\ell) L_{k+\ell} .
$$
Note that this is a representation of the Lie subalgebra
$\Vir_+\subset\Vir$ of the Virasoro algebra spanned by $L_k$, $k\ge-1$; it
is isomorphic to the the Lie algebra of polynomial coefficient vector
fields on the line, by the map $L_k\mapsto-x^{k+1}\,\p/\p x$.

If $X$ is a point, the generating function $Z=Z^X$ equals
\begin{equation}
Z = \exp \biggl( \sum_{g\ge 0} \hbar^{g-1} \sum_{n\ge 0} \frac{1}{n!}
\sum_{k_1\dots k_n} t_{k_1} \dots t_{k_n} \< \tau_{k_1} \dots \tau_{k_n}
\>_g \biggr) .
\end{equation}
where
$$
\< \tau_{k_1} \dots \tau_{k_n} \>_g = \int_{\Mbar_{g,n}} \psi_1^{k_1} \dots
\psi_n^{k_n} .
$$
A conjecture of Witten, proved by Kontsevich \cite{K}, asserts the
annihilation of $Z$ by the operators
$$
L_k = \begin{cases}
\displaystyle \sum_{m=1}^\infty \bigl( t_m - \delta_{m1} \bigr) \p_{m-1} +
\frac{1}{2\hbar} t_0^2 , & k=-1 , \\
\displaystyle \sum_{m=0}^\infty (m+\half) \bigl( t_m - \delta_{m1} \bigr)
\p_m + \frac{1}{16} , & k=0 , \\
\displaystyle \sum_{m=0}^\infty
\tfrac{\Gamma(k+m+\frac{3}{2})}{\Gamma(m+\frac{1}{2})} \bigl( t_m -
\delta_{m1} \bigr) \p_{m+k} + \frac{\hbar}{2} \sum_{m=0}^{k-1} (-1)^{m+1}
\tfrac{\Gamma(k-m+\frac{1}{2})}{\Gamma(-m-\frac{1}{2})} \p_m
\p_{k-m-1} , & k>0 .
\end{cases}
$$
The Virasoro conjecture for arbitrary smooth projective varieties $X$
may be viewed as a generalization of this conjecture of Witten.

The Virasoro conjecture differs in some respects from Witten's conjecture:
whereas Witten's conjecture determines all intersection numbers
$\<\tau_{k_1}\dots\tau_{k_n}\>_g$ in terms of the basic one
$\<\tau_0^3\>_0=1$, the Virasoro conjecture does not appear to suffice to
determine the descendent Gromov-Witten invariants of positive-dimensional
$X$. Furthermore, although Gromov-Witten invariants may be defined for any
compact symplectic manifold $X$, the Virasoro conjecture depends on the
Hodge decomposition of $H^\bull(X,\C)$, and thus does not appear to
generalize beyond K\"ahler manifolds. Furthermore, the action of $\Vir_+$
only extends to an action of $\Vir$ when $X$ is even dimensional (in which
case it has central charge the Euler characteristic $\chi(X)$ of $X$).

In \cite{EHX1}, the authors outline a proof of the Virasoro conjecture
in genus $0$ using the genus $0$ topological recursion
relation. Throughout this paper, we will assume that the Virasoro
conjecture holds in genus $0$.

There is a natural conjecture lying between those of Witten and of
Eguchi-Hori-Xiong. The Virasoro conjecture for $X$ implies that the
generating function for degree $0$ descendent Gromov-Witten invariants
$$
Z^X_0 = \exp \biggl( \sum_{g\ge 0} \hbar^{g-1} F^X_{g,\beta=0} \biggr)
$$
is also annihilated by the representation $\rho^X$. It is the
implications of this \punctual Virasoro conjecture, in genus $g>0$,
that we investigate here.

\subsection{Chern classes of the Hodge bundle}
Let $g>1$. The zero degree Gromov-Witten invariants which we study are
obtained by integrating against a cycle of dimension $(r-3)(1-g)+n$ in the
moduli stack $\Mbar_{g,n}(X,0)$ of stable maps, called the virtual
fundamental class. Axiom~IV for Gromov-Witten invariants of Behrend
\cite{B} states that this virtual fundamental class is the flat pullback by
the map $\Mbar_{g,n}(X,0)\to\Mbar_{g,0}(X,0)$ of the virtual fundamental
class of $\Mbar_{g,0}(X,0)$; in particular, it vanishes if $r>3$. Thus,
there are only three cases to be considered: $X$ is a curve, a surface or a
threefold. (The case where $X$ is zero-dimensional is precisely Witten's
conjecture.)

Let $\Mbar_{g,n+1}/\Mbar_{g,n}$ be the universal curve $n$-pointed curve of
genus $g$, let $\pi_{g,n}$ be the projection from $\Mbar_{g,n+1}$ to
$\Mbar_{g,n}$, and let $\om_{g,n}=\om_{\Mbar_{g,n+1}/\Mbar_{g,n}}$ be the
relative dualizing sheaf. The Hodge bundle on $\Mbar_{g,n}$ is the vector
bundle $\EE_{g,n}=\pi_{g,n*}\om_{g,n}$. Note that
\begin{equation} \label{gn}
\pi_{g,n}^*\EE_{g,n} \cong \EE_{g,n+1} .
\end{equation}

Fix a genus $g$, let $\lambda_i$ be the $i$th Chern class $c_i(\EE)$ of
$\EE=\EE_g$, and let $c_t(\EE)$ be the total Chern class of $\EE$:
$$
c_t(\EE) = \sum_{i=0}^g t^i \lambda_i .
$$
The omission of the number of marked points $n$ from the notation for
$\lambda_i$ is justified, since $\pi_{g,n}^*c_i(\EE_{g,n}) =
c_i(\EE_{g,n+1})$ by \eqref{gn}.

Mumford \cite{Mumford} has proved the relation $c_t(\EE) c_{-t}(\EE) = 1$.
Extracting the coefficients of $t^{2g}$ and $t^{2g-1}$, we see in
particular that
\begin{equation}
\begin{cases}
\lambda_g^2=0 , & \\
\lambda_{g-1}^2 = 2\lambda_g\lambda_{g-2} . &
\end{cases}
\end{equation}

Since $\Mbar_{g,n}(X,0)$ is isomorphic to $X\times\Mbar_{g,n}$, it has
dimension $r+(3g-3+n)$, and the virtual fundamental class equals
$e(T_X\boxtimes\EE_g^\Dual)\cap[X\times\Mbar_{g,n}]$, where
$e(T_X\boxtimes\EE_g^\Dual)=c_{rg}(T_X\boxtimes\EE_g^\Dual)$ is the Euler
class, or top Chern class, of the rank $rg$ bundle
$T_X\boxtimes\EE_g^\Dual$. It is for this reason that the \punctual
Virasoro conjecture involves intersection numbers of $\psi$ and $\lambda$
classes.

\subsection{The \punctual Virasoro conjecture for curves}
We now turn to a summary of our results. The \punctual Virasoro conjecture
for $\P^1$ implies that if $2g-3+n=k_1+\dots+k_n$, then%
\footnote{C.~Faber helped us to find this explicit expression.}
\begin{equation} \label{lamg}
\int_{\Mbar_{g,n}} \psi_1^{k_1} \dots \psi_n^{k_n} \lambda_g =
\binom{2g+n-3}{k_1,\ldots, k_n} b_g ,
\end{equation}
where
$$
b_g = \begin{cases} 1 , & g=0 , \\ \displaystyle \int_{\Mbar_{g,1}}
\psi_1^{2g-2} \lambda_g , & g>0 .
\end{cases}
$$
Equation \eqref{lamg} generalizes the well-known multinomial formula for
the $\psi$ integrals in genus $0$. The constants $b_g$, $g>0$, are
calculated in \cite{FP} using algebro-geometric techniques:
\begin{equation} \label{ground1} 
\sum_{g=0}^\infty b_g t^{2g} = \frac{t/2}{\sin(t/2)} ,
\end{equation}
or, in terms of Bernoulli numbers,
$$
b_g = \frac{2^{2g-1}-1}{2^{2g-1}} \frac{|B_{2g}|}{(2g)!} .
$$
They do not appear to be constrained by the Virasoro conjecture, \punctual
or otherwise.

The \punctual Virasoro conjecture for a curve is equivalent to \eqref{lamg}
together with an explicit recursion relation (see Section \ref{P1}) for the
intersection numbers
$$
\int_{\Mbar_{g,n}} \psi_1^{k_1} \ldots \psi_n^{k_n} \lambda_{g-1} ,
$$
for which we do not know a closed solution. In particular, in combination
with (\ref{ground1}), the conjecture implies that
\begin{equation} \label{cg}
c_g = \int_{\Mbar_{g,1}} \psi_1^{2g-1} \lambda_{g-1} = \Bigl(
\sum_{k=1}^{2g-1} \frac{1}{k} \Bigr) b_g - \frac12 \sum_{g=g_1+g_2}
\frac{(2g_1-1)! (2g_2-1)!}{(2g-1)!} b_{g_1} b_{g_2} .
\end{equation}

All integrals of $\psi$ and $\lambda$ classes over $\Mbar_{g,n}$ may in
principle be calculated by an algorithm of Faber \cite{F2}, which he has
implemented in \texttt{Maple}. This algorithm uses Mumford's
Grothendiek-Riemann-Roch formulas to replace factors of $\lambda_i$ in the
integrand by combinations of boundary divisor and $\psi$ classes. The
resulting integrals may then be reduced to pure $\psi$ integrals, which may
be calculated by Witten's conjectures. Unfortunately, it appears to be
impractical to prove the \punctual Virasoro conjecture using this
algorithm.

We have verified \eqref{cg} up to genus $5$ using Faber's program,
obtaining the following results:
$$\begin{tabular}{|C||C|C|} \hline
g & b_g & c_g \\ \hline
1 & 1/24 & 1/24 \\
2 & 7/5760 & 1/480 \\
3 & 31/967680 & 41/580608 \\
4 & 127/154828800 & 13/6220800 \\
5 & 73/3503554560 & 21481/367873228800 \\ \hline
\end{tabular}$$

\medskip

\subsection{The \punctual Virasoro conjecture for surfaces}
The \punctual Virasoro conjecture for $\P^2$ implies that if
$g-1+n=k_1+\dots+k_n$ and $k_i>0$,
\begin{equation} \label{lamgg}
\int_{\Mbar_{g,n}} \psi_1^{k_1} \dots \psi_n^{k_n} \lambda_g \lambda_{g-1}
= \frac{(2g+n-3)! (2g-1)!!}{(2g-1)! (2k_1-1)!!  \dots (2k_n-1)!!}
\int_{\Mbar_{g,1}} \psi_1^{g-1} \lambda_g \lambda_{g-1} .
\end{equation}
The constant
\begin{equation} \label{ground2} 
\int_{\Mbar_{g,1}} \psi_1^{2g-2} \lambda_g \lambda_{g-1} =
\frac{1}{2^{2g-1}(2g-1)!!} \frac{|B_{2g}|}{2g}
\end{equation}
has been calculated by Faber \cite{F3}, who shows that it follows from
Witten's conjecture.

Remarkably, \eqref{lamgg} is part of a deep conjecture of Faber \cite{F1}
concerning the so-called tautological Chow ring $R^\bull(\CM_g)$ generated
over $\Q$ by the classes $\kappa_i$. Combining results of Looijenga
\cite{L} and Faber \cite{F3}, we know that if $i_1+\dots+i_m=g-2$, then
$$
\kappa_{i_1} \dots \kappa_{i_m} = \frac{\displaystyle\int_{\Mbar_g}
\kappa_{i_1} \dots \kappa_{i_m} \lambda_g\lambda_{g-1}}
{\displaystyle\int_{\Mbar_g} \kappa_{g-2} \lambda_g\lambda_{g-1}}
\kappa_{g-2} \in A^{g-2}(\CM_g)_{\Q} .
$$
By the formula
$$
\int_{\Mbar_{g,n}} \psi_1^{k_1+1} \dots \psi_n^{k_n+1} \lambda_g
\lambda_{g-1} = \sum_{\sigma\in S_n} \int_{\Mbar_g} \kappa_\sigma \lambda_g
\lambda_{g-1} ,
$$
where $\kappa_\sigma$ is the product of $\kappa_{|\CO|}$, one for each
cycle $\CO$ of $\sigma$, and $|\CO|=\sum_{i\in\CO}k_i$, the integrals in
this expression may be calculated from \eqref{lamgg}. In \cite{F1}, Faber
proves \eqref{lamgg} for genus $g\le 15$, and conjectures that it holds in
all genera.

Assuming \eqref{lamgg}, the \punctual Virasoro conjecture for surfaces is
equivalent to a complicated recursion relation for the intersection numbers
$$
\int_{\Mbar_{g,n}} \psi_1^{k_1} \ldots \psi_n^{k_n} \lambda_g \lambda_{g-2} ,
$$
for which we do not know a closed solution.

\subsection{The \punctual Virasoro conjecture for threefolds}
The \punctual Virasoro conjecture for threefolds is already implied by the
string and dilaton equations: the constants in this case are the
intersection numbers
\begin{equation} \label{cycont}
\int_{\Mbar_g} \lambda_{g-1}^3 = \frac{1}{(2g-2)!} \frac{|B_{2g-2}|}{2g-2}
\frac{|B_{2g}|}{2g} ,
\end{equation}
whose values were conjectured by Faber \cite{F1}, and are calculated in
\cite{FP}.

\section{The Euler characteristic of the obstruction bundle}

The moduli space of degree $0$ maps to $X$ has a very simple form:
\begin{equation} \label{degzero} 
\Mbar_{g,n}(X,0) = X \times \Mbar_{g,n}.
\end{equation}

The virtual fundamental class $[\Mbar_{g,n}(X, 0)]^{\text{vir}}$ is
equal to $e(T_X\boxtimes\EE^\Dual)\cap[X\times\Mbar_{g,n}]$ via the
isomorphism \eqref{degzero}, where
$e(T_X\boxtimes\EE^\Dual)=c_{rg}(T_X\boxtimes\EE^\Dual)$ denotes the
Euler class, or top Chern class, of the vector bundle
$T_X\boxtimes\EE^\Dual$. Let $\{\gamma_a\}$ be a basis of $H^*(X,\Q)$.
Let $\psi_i$ denote the first Chern class of the $i$th cotangent line
bundle on the moduli space of maps.  The degree $0$ gravitational
descendents of $X$ are the integrals:
$$
\< \tau_{k_1}(\gamma_{a_1}) \dots \tau_{k_n}(\gamma_{a_n}) \>_{g,0}^X =
\int_{X\times\Mbar_{g,n}} \gamma_{a_1} \dots \gamma_{a_n} \psi_1^{k_1}
\dots \psi_n^{k_n} \cup e(T_X\boxtimes\EE^\Dual) .
$$
These descendents involve the cohomology ring of $X$ and the integrals over
$\Mbar_{g,n}$ of the $\psi$ and $\lambda$ classes. The \punctual Virasoro
conjectures imply relations among the latter set of integrals on
$\Mbar_{g,n}$.

In this section, we calculate the Euler class $e(\CE)$ of the
obstruction bundle $\CE=T_X\boxtimes\EE^\Dual$ on the moduli space
$\Mbar_{g,n}(X,0)\cong X\times\Mbar_{g,n}$ of degree $0$ stable maps
in the three cases in which $X$ is a curve, a surface and a threefold.

The case $g=1$ is exceptional, since the Euler class $e(\CE)$ may be
nonzero for $X$ of any dimension: it is easily seen that
$$
e(\CE) = c_r(X) - c_{r-1}(X) \lambda_1 .
$$

\subsection{$X$ a curve}
Observe that if $L$ is a line bundle,
$$
c(L\boxtimes\EE^\Dual) = \sum_{i=0}^g c_1(L)^i c_{g-i}(\EE^\Dual) 
= \sum_{i=0}^g (-1)^{g-i} c_1(L)^i \lambda_{g-i} .
$$
If $X$ is a curve, we may set $L=T_X$, so that $c_1(L)=c_1(X)$. Since
$c_1(L)^i$ vanishes if $i>1$, we conclude that
$$
(-1)^g e(\CE) = \lambda_g - c_1(X) \lambda_{g-1} .
$$

\subsection{$X$ a surface}
In this case, by the splitting principle, we may suppose that $T_X\cong
L_1\oplus L_2$ is the sum of two line bundles. We see that
\begin{align*}
e(\CE) &= e(L_1\boxtimes\EE^\Dual)e(L_2\boxtimes\EE^\Dual) \\ &= \bigl(
\lambda_g - c_1(L_1) \lambda_{g-1} + c_1(L_1)^2 \lambda_{g-2} \bigr) \bigl(
\lambda_g - c_1(L_2) \lambda_{g-1} + c_1(L_2)^2 \lambda_{g-2} \bigr) \\ &=
\lambda_g^2 - c_1(X) \lambda_g\lambda_{g-1} + \bigl( c_1(X)^2 - 2 c_2(X)
\bigr) \lambda_g \lambda_{g-2} + c_2(X) \lambda_{g-1}^2 .
\end{align*}
Applying Mumford's relations $\lambda_g^2=0$ and
$\lambda_{g-1}^2=2\lambda_g\lambda_{g-2}$, we see that
$$
e(\CE) = - c_1(X) \lambda_g\lambda_{g-1} + c_1(X)^2 \lambda_g\lambda_{g-2}
.
$$

\subsection{$X$ a threefold}

By the splitting principle, we may suppose that the tangent bundle
$T_X\cong L_1\oplus L_2\oplus L_3$ is the sum of three line bundles. We see
that
\begin{align*}
e(\CE) &=
e(L_1\boxtimes\EE^\Dual)e(L_2\boxtimes\EE^\Dual)e(L_3\boxtimes\EE^\Dual) \\
&= (-1)^g \bigl( \lambda_g - c_1(L_1) \lambda_{g-1} +
c_1(L_1)^2\lambda_{g-2} - c_1(L_1)^3\lambda_{g-3} \bigr) \\ & \quad \bigl(
\lambda_g - c_1(L_2) \lambda_{g-1} + c_1(L_1)^2\lambda_{g-2} -
c_1(L_2)^3\lambda_{g-3} \bigr) \\ & \quad \bigl( \lambda_g - c_1(L_3)
\lambda_{g-1} + c_1(L_3)^2\lambda_{g-2} - c_1(L_3)^3\lambda_{g-3} \bigr) .
\end{align*}
Since $\lambda_g^2=0$ and
$\lambda_g\lambda_{g-1}^2=2\lambda_g^2\lambda_{g-2}=0$, many terms of the
expansion of this product drop out, and we see that
\begin{align*}
(-1)^g e(\CE) &= - \sum_{i\ne j} c_1(L_i) c_1(L_j)^2
\lambda_g\lambda_{g-1}\lambda_{g-1} - c_3(X) \lambda_{g-1}^3 \\
&= \bigl( 3c_3(X) - c_2(X)c_1(X) \bigr) \lambda_g\lambda_{g-1}\lambda_{g-2}
- c_3(X) \lambda_{g-1}^3 \\ 
&= \tfrac{1}{2} \bigl( c_3(X) - c_2(X)c_1(X) \bigr) \lambda_{g-1}^3 .
\end{align*}

\section{The Virasoro conjecture} \label{Virasoro}

Let $X$ be a smooth projective variety of dimension $r$, and let $\gamma_a$
be a basis for $H^\bull(X,\C)$; we suppose that the cohomology classes are
homogeneous with respect to the Hodge decomposition, so that there exist
integers $p_a$ and $q_a$ such that $\gamma_a\in H^{p_a,q_a}(X)$. Let
$b_a=p_a+(1-r)/2$.

In the following formulas, we use the Einstein summation convention over
indices $a$ and $b$, making use of the non-degenerate inner product
$$
\eta_{ab} = \int_X \gamma_a \cup \gamma_b
$$
and its inverse $\eta^{ab}$ to raise and lower indices as needed. Let
$C^b_a$ be the matrix of the first Chern class of $X$:
$$
C^b_a \gamma_b = c_1(X) \cup \gamma_a .
$$

Introduce an affine space with coordinates $\{t^a_k\mid k\ge0\}$, called
the large phase space; the full Gromov-Witten potential \eqref{GW} is a
formal function on this space. Let $\p_{a,k}=\p/\p t^a_k$, and let
$\tilde{t}^a_k = t^a_k-\delta_{a0}\delta_{k1}$. Let
$$
[x]^k_i = e_{k+1-i}(x,x+1,\dots,x+k) ,
$$
where $e_k$ is the $k$th elementary symmetric function of its arguments;
thus,
$$
\sum_{i=0}^{k+1} [x]^k_i t^i = (t+x)(t+x+1)\dots(t+x+k) .
$$

Following Eguchi, Hori and Xiong \cite{EHX1,EHX2}, we introduce
differential operators $L_k$, $k\ge-1$, by the formulas
\begin{align*}
L_k &= \sum_{m=0}^\infty \sum_{i=0}^{k+1} \Bigl( [b_a\!+\!m]^k_i (C^i)^b_a
\tilde{t}^a_m \p_{b,m+k-i} + \frac{\hbar}{2} (-1)^{m+1}
[b_a\!-\!m\!-\!1]^k_i (C^i)^{ab} \p_{a,m} \p_{b,k-m-i-1} \Bigr) \\ & \quad
+ \frac{1}{2\hbar} (C^{k+1})_{ab} t^a_0 t^b_0 + \frac{\delta_{k0}}{48}
\int_X \bigl( (3-r)c_r(X)-2c_1(X)c_{r-1}(X) \bigr) ,
\end{align*}
where it is understand that $\tilde{t}^a_m$ and $\p_{a,m}$ vanish if
$m<0$. Note that the conjecture of Eguchi-Hori-Xiong is for projective
varieties such that $p_a=q_a$ for all $a$; the extension of their
conjecture to general smooth projective varieties is due to Katz (private
communication, March 1997).


\section{The \punctual Virasoro conjecture for a curve} \label{P1}

Let $X$ be a smooth projective curve of genus $\gamma$. Choose dual bases
$(e^1,\dots,e^g)$ and $(f^1,\dots,f^g)$ of $H^{0,1}(X)$ and
$H^{1,0}(X)$. Denote by $t_k$, $k\ge0$, the coordinates on the large phase
space dual to the descendents $\tau_k(1)$ of $1\in H^0(X)$, by
$\Alpha_k=(\alpha_k^1,\dots,\alpha_k^g)$ and
$\Beta_k=(\beta_k^1,\dots,\beta_k^g)$ the coordinates dual to the
descendents $\tau_k(e^i)$ and $\tau_k(f^i)$, and by $s_k$ the coordinates
dual to the descendents $\tau_k(\om)$ of the class $\om\in H^2(X)$
Poincar\'e dual to a point.

If $\alpha$ is a cohomology class on $\Mbar_g$ (and hence, by pullback, on
the moduli spaces $\Mbar_{g,n}$), introduce the generating functions
\begin{equation} \label{psilambda}
\<\<\tau_{k_1}\dots\tau_{k_n}\mid\alpha\>\>_g = \sum_{N=0}^\infty
\frac{1}{N!}  \sum_{l_1\dots l_N} t_{l_1} \dots t_{l_N}
\int_{\Mbar_{g,n+N}} \psi^{k_1}\dots\psi^{k_n}\psi^{l_1}\dots\psi^{l_N}
\alpha .
\end{equation}

\begin{theorem}
We have
$$
\frac{L_kZ^X_0}{Z^X_0} = \sum_{g=0}^\infty \hbar^{g-1} (-1)^g \biggl(
(2\gamma-2) x^k_g(\t) + \sum_{\ell=0}^\infty \Bigl( s_\ell +
\sum_{m=0}^\infty \Alpha_\ell \* \Beta_m \p_{t_m} \Bigr) y^k_{g,\ell}(\t)
\biggr) ,
$$
where
\begin{align*}
x^k_g(\t) & = -[1]^k_0 \<\<\tau_{k+1}\mid\lambda_{g-1}\>\>_g +
\sum_{m=0}^\infty t_m [m]^k_0 \<\<\tau_{k+m}\mid\lambda_{g-1}\>\>_g \\ &
\quad + [1]^k_1 \<\<\tau_k\mid\lambda_g\>\>_g - \sum_{m=0}^\infty t_m [m]^k_1
\<\<\tau_{k+m-1}\mid\lambda_g\>\>_g \\ & \quad - \frac12 \sum_{m=0}^{k-2}
\sum_{g=g_1+g_2} (-1)^{m+1} [-m\!-\!1]^k_1 \<\<\tau_m\mid\lambda_{g_1}\>\>_{g_1}
\<\<\tau_{k-m-2}\mid\lambda_{g_2}\>\>_{g_2} , \\
y^k_{g,\ell}(\t) & = - [1]^k_0 \<\<\tau_{k+1}\tau_\ell\mid\lambda_g\>\>_g +
\sum_{m=1}^\infty [m]^k_0 t_m \<\<\tau_{k+m}\tau_\ell\mid\lambda_g\>\>_g +
[\ell\!+\!1]^k_0 \<\<\tau_{k+\ell}\mid\lambda_g\>\>_g
\end{align*}
\end{theorem}
\begin{proof}
For $k>0$, $L_k$ is given by the formula
\begin{multline*}
L_k = - [1]^k_0 \, \p_{t_{k+1}} + \sum_{m=0}^\infty \Bigl( [m]^k_0 \bigl(
t_m \p_{t_{k+m}} + \Alpha_m \* \p_{\Alpha_{k+m}} \bigr) + [m\!+\!1]^k_0
\bigl( s_m \p_{s_{k+m}} + \Beta_m \* \p_{\Beta_{k+m}} \bigr) \Bigr) \\ +
(2-2\gamma) \biggl( - [1]^k_{1} \, \p_{s_k} + \sum_{m=0}^\infty [m]^k_1 t_m
\p_{s_{k+m-1}} + \frac{\hbar}{2} \sum_{m=0}^{k-2} (-1)^{m+1} [-m\!-\!1]^k_1
\p_{s_m} \p_{s_{k-m-2}} \biggr) .
\end{multline*}
Using the notation \eqref{psilambda}, $Z^X_0$ is given by the formula
\begin{multline*}
Z^X_0 = \exp \Biggl( (2\gamma-2) \sum_{g=1}^\infty (-1)^g \hbar^{g-1}
\<\<~\mid\lambda_{g-1}\>\>_g \\ + \sum_{g=0}^\infty (-1)^g \hbar^{g-1}
\biggl( \sum_{m=0}^\infty s_m \<\<\tau_m\mid\lambda_g\>\>_g +
\sum_{\ell,m=0}^\infty \Alpha_\ell \* \Beta_m
\<\<\tau_\ell\tau_m\mid\lambda_g\>\>_g \biggr) \Biggr) .
\end{multline*}
The theorem follows on combining these formulas.
\end{proof}

\begin{corollary}
The \punctual Virasoro conjecture for algebraic curves is equivalent to the
vanishing of the expressions $x^k_g(\t)$ and $y^k_{g,\ell}(\t)$. In
particular, if the \punctual Virasoro conjecture holds for $\P^1$, then it
holds for all curves.
\end{corollary}

The vanishing of $y^k_{g,\ell}(\t)$ for $k\ge1$ and $\ell\ge0$ is
equivalent to the formula \eqref{lamg} for the generating function
$\<\<~\mid\lambda_g\>\>_g$. To see this, note that
\begin{multline} \label{recurse:P1}
\frac{1}{(k+1)!} \p_{t_{k_1}} \dots \p_{t_{k_n}} y^k_{g,k_0}(\mathbf{0}) =
- \<\tau_{k+1}\tau_{k_0}\dots\tau_{k_n}\mid\lambda_g\>_g \\ +
\binom{k_0+k+1}{k_0}
\<\tau_{k_0+k}\tau_{k_1}\dots\tau_{k_n}\mid\lambda_g\>_g + \sum_{i=1}^n
\binom{k_i+k}{k_i-1}
\<\tau_{k_0}\dots\tau_{k_i+k}\dots\tau_{k_n}\mid\lambda_g\>_g ,
\end{multline}
where it is understood that $\binom{a}{-1}=0$ for $a$ a natural number.
\begin{theorem} \label{details}
The recursion given by the vanishing of \eqref{recurse:P1} has the unique
solution
$$
\<\tau_{k_1}\dots\tau_{k_n}\mid\lambda_g\>_g =
\begin{cases}
\displaystyle \binom{n-3}{k_1,\ldots,k_n} \<\tau_0^3\>_0 , & g=0 , \\[10pt]
\displaystyle \binom{2g+n-3}{k_1,\ldots,k_n} \<\tau_{2g-2}\mid\lambda_g\> , & g>0 .
\end{cases}$$
\end{theorem}
\begin{proof}
We prove the theorem by induction on $n$; in the cases $n=3$ for $g=0$ and
$n=1$ for $g>0$, the formula is a tautology. Thus, we must prove that
\begin{multline*}
\binom{2g+n-1}{k_0,\dots,k_n,k+1} = \binom{k_0+k+1}{k_0}
\binom{2g+n-2}{k_0+k,k_1,\dots,k_n} \\ + \sum_{i=1}^n \binom{k_i+k}{k_i-1}
\binom{2g+n-2}{k_0,\dots,k_i+k,\dots,k_n} .
\end{multline*}
This follows from the equation
$$
2g+n-1 = (k_0+k+1) + \sum_{i=1}^n k_i ,
$$
on multiplication of both sides by $\dfrac{(2g+n-2)!}{k_0!\dots k_n!
(k+1)!}$.
\end{proof}

In particular, the well-known formulas for the intersection numbers
$\<\tau_{k_1}\dots\tau_{k_n}\>_0$ are seen to be special cases of the
conjectured formulas for the intersection numbers
$\<\tau_{k_1}\dots\tau_{k_n}\mid\lambda_g\>_g$.

The same technique applied to $x^k_g(\t)$ leads to a recursion for the
intersection numbers $\<\tau_{k_1}\dots\tau_{k_n}\mid\lambda_{g-1}\>_g$
which expresses them in terms of the numbers $b_h$, $h\le g$. We will only
discuss the simplest case $n=1$. Taking the relation
$x^{2g-2}_g(\mathbf{0})=0$, we obtain the formula
$$
(2g-1)! c_g = s(2g,2) b_g - \frac12 \sum_{g=g_1+g_2} (2g_1-1)! (2g_2-1)!
b_{g_1} b_{g_2} ,
$$
where $s(2g,2)$ is the Stirling number of the first kind
$$
s(2g,2) = [1]^{2g-2}_1 = (2g-1)! \sum_{k=1}^{2g-1} \frac{1}{k} .
$$
We have not been able to find a solution of this recursion, or its
generalizations to larger $n$, in closed form.

\raggedbottom

\section{The \punctual Virasoro conjecture for a surface}

The discussion of the \punctual Virasoro conjecture for a surface runs
along the same lines as for a curve, although the details are a little more
complicated. To simplify notation, we restrict attention to
simply-connected surfaces.

Let $X$ be a smooth simply-connected projective surface. Choose dual bases
$(e_1,\dots,e_p)$ and $(f_1,\dots,f_p)$ of $H^{0,2}(X)$ and $H^{2,0}(X)$,
and a basis $\om_i$, $1\le i\le\ell$ of $H^{1,1}(X)$. Denote by $t_k$,
$\s_k=(s^1_k,\dots,s^d_k)$, $r_k$, $\a_k=(a^1_k,\dots,a^p_k)$ and
$\b_k=(b^1_k,\dots,b^p_k)$ the coordinates on the large phase space dual
respectively to the descendents of $1\in H^0(X,\Z)$, of $\om_i$, $1\le i\le
d$, of the class in $H^4(X,\Z)$ Poincar\'e dual to a point, and of $e_i$
and $f_i$, $1\le i\le p$, respectively. Let $\c=(c_1,\dots,c_d)$ be the
vector in the vector space dual to $H^{1,1}(X)$ such that $c_1(X)=\c\*\om$.

\begin{theorem}
We have
$$
\frac{L_kZ^X_0}{Z^X_0} = \sum_{g=1}^\infty \hbar^{g-1} (-1)^g \biggl(
|\c|^2 x^k_g(\t) - \sum_{\ell=0}^\infty \c\*\s_\ell \, y^k_{g,\ell}(\t)
\biggr) + \frac{1}{\hbar} w(r,\a,\s,\b,t)
$$
where
\begin{align*}
x^k_g(\t) & = - [\half]^k_0 \<\<\tau_{k+1}\mid\lambda_g\lambda_{g-2}\>\>_g -
\sum_{m=0}^\infty t_m [m\!-\!\half]^k_0
\<\<\tau_{k+m}\mid\lambda_g\lambda_{g-2}\>\>_g \\
& \quad + \sum_{m=0}^{k-1} (-1)^{m+1} \Bigl( [-m\!-\!\tfrac32]^k_0
\<\<\tau_m\>\>_0 \<\<\tau_{k-m-1}\mid\lambda_g\lambda_{g-2}\>\>_g \\
& \qquad + \half [-m\!-\!\half]^k_0 \sum_{g=g_1+g_2}
\<\<\tau_m\mid\lambda_{g_1}\lambda_{g_1-2}\>\>_{g_1}
\<\<\tau_{k-m-1}\mid\lambda_{g_2}\lambda_{g_2-2}\>\>_{g_2} \Bigr) \\
& \quad + [\half]^k_1 \<\<\tau_k\mid\lambda_g\lambda_{g-1}\>\>_g +
\sum_{m=0}^\infty t_m [m\!-\!\half]^k_1
\<\<\tau_{k+m-1}\mid\lambda_g\lambda_{g-1}\>\>_g \\ & \quad +
\sum_{m=0}^{k-2} (-1)^{m+1} [-m\!-\!\tfrac32]^k_1 \<\<\tau_m\>\>_0
\<\<\tau_{k-m-2}\mid\lambda_g\lambda_{g-1}\>\>_g , \\
y^k_{g,\ell} & = - [\half]^k_0
\<\<\tau_{k+1}\tau_\ell\mid\lambda_g\lambda_{g-1}\>\>_g + \sum_{m=0}^\infty
t_m [m\!-\!\half]^k_0 \<\<\tau_{k+m}\tau_\ell\mid\lambda_g\lambda_{g-1}\>\>_g \\
& \quad + [\ell\!+\!\half]^k_0 \<\<\tau_{k+\ell}\mid\lambda_g\lambda_{g-1}\>\>_g \\
& \quad + \sum_{m=0}^{k-1} (-1)^{m+1} \Bigl( [-m\!-\!\tfrac32]^k_0
\<\<\tau_m\>\>_0 \<\<\tau_{k-m-1}\tau_\ell\mid\lambda_g\lambda_{g-1}\>\>_g \\
& \qquad + [-m\!-\!\half]^k_0 \<\<\tau_m\tau_\ell\>\>_0
\<\<\tau_{k-m-1}\mid\lambda_g\lambda_{g-1}\>\>_g \Bigr)
\end{align*}
\end{theorem}
We have omitted the explicit expression for $w(r,\a,\s,\b,t)$, because of
its greater complexity, because it differs in nature from the higher genus
coefficients, and because in any case we are assuming that it vanishes by
the genus $0$ Virasoro conjecture.
\begin{proof}
For $k>0$, $L_k$ is given by the formula
\begin{align*}
L_k &= - [\half]^k_0 \p_{t_{k+1}} + \sum_{m=0}^\infty \Bigl(
[m\!-\!\half]^k_0 \bigl( t_m \p_{t_{k+m}} + \b_m \* \p_{\b_{k+m}} \bigr) +
[m\!+\!\half]^k_0 \s_m \* \p_{\s_{k+m}} \\
{} & \qquad\qquad\qquad\qquad + [m\!+\!\tfrac{3}{2}]^k_0 \bigl( r_m
\p_{r_{k+m}} + \a_m \* \p_{\a_{k+m}} \bigr) \Bigr) \\
{} & \qquad\qquad + \hbar \sum_{m=0}^{k-1} (-1)^{m+1} \Bigl(
[-m\!-\!\tfrac{3}{2}]^k_0 \p_{r_m} \p_{t_{k-m\!-\!1}} + \half
[-m\!-\!\half]^k_0 \p_{\s_m} \* \p_{\s_{k-m\!-\!1}} \Bigr) \\
{} & + \c \* \biggl( - [\half]^k_1 \p_{\s_k} + \sum_{m=0}^\infty \bigl(
[m\!-\!\half]^k_1 t_m \p_{\s_{k+m\!-\!1}} + [m\!+\!\half]^k_1 \s_m
\p_{r_{k+m\!-\!1}} \bigr) \\
{} & \qquad\qquad + \hbar \sum_{m=0}^{k-2} (-1)^{m+1}
[-m\!-\!\tfrac{3}{2}]^k_1 \p_{r_m} \p_{\s_{k-m\!-\!2}} \biggr) \\
{} & + |\c|^2 \biggl( - [\half]^k_2 \p_{r_{k-1}} + \sum_{m=0}^\infty
[m\!-\!\half]^k_2 t_m \p_{r_{k+m-2}} \\
{} & \qquad\qquad + \frac{\hbar}{2} \sum_{m=0}^{k-3} (-1)^{m+1}
[-m\!-\!\tfrac{3}{2}]^k_2 \p_{r_m} \p_{r_{k-m-3}} +
\frac{\delta_{k1}}{2\hbar} t_0^2 \biggr)
\end{align*}
The generating function $Z^X_0$ is given by the formula
\begin{multline*}
Z^X_0 = \exp \Bigl( \frac{1}{\hbar} \sum_{m=0}^\infty r_m \<\<\tau_m\>\>_0
+ \frac{1}{\hbar} \sum_{\ell,m} \bigl( \half \s_\ell\*\s_m + \a_\ell\*\b_m \bigr)
\<\<\tau_\ell\tau_m\>\>_0 \\ + |\c|^2 \sum_{g=1}^\infty \hbar^{g-1}
\<\<~\mid\lambda_g\lambda_{g-2}\>\>_g - \sum_{g=1}^\infty \hbar^{g-1}
\sum_{m=0}^\infty \c\*\s_m \<\<\tau_m\mid\lambda_g\lambda_{g-1}\>\>_g
\Bigr)
\end{multline*}
The theorem follows on combining these formulas.
\end{proof}

\begin{corollary}
The \punctual Virasoro conjecture for surfaces is equivalent to the
vanishing of the expressions $x^k_g(\t)$ and $y^k_{g,\ell}(\t)$. In
particular, if the \punctual Virasoro conjecture holds for $\P^2$, then it
holds for all simply connected surfaces.
\end{corollary}

The vanishing of $y^k_{g,\ell}(\t)$ for $k\ge0$ and $\ell>0$ implies
formula \eqref{lamgg} for the generating function
$\<\<~\mid\lambda_g\lambda_{g-1}\>\>_g$. To see this, note that if
$k_1,\dots,k_n>0$,
\begin{align} \notag
\frac{1}{[\half]^k_0} \p_{t_{k_1}} \dots \p_{t_{k_n}} y_{g,k_0}(\mathbf{0}) &=
- \<\tau_{k+1}\tau_{k_0}\dots\tau_{k_n}\mid\lambda_g\lambda_{g-1}\>_g
\\ & \quad + \frac{(2k+2k_0+1)!!}{(2k+1)!!(2k_0-1)!!}
\<\tau_{k_0+k}\tau_{k_1}\dots\tau_{k_n}\mid\lambda_g\lambda_{g-1}\>_g
\label{recurse:P2} \\ & \quad +
\sum_{i=1}^n \frac{(2k+2k_i-1)!!}{(2k+1)!!(2k_i-3)!!}
\<\tau_{k_0}\dots\tau_{k_i+k}\dots\tau_{k_n}\mid\lambda_g\lambda_{g-1}\>_g
\notag
\end{align}
The proof of the following theorem is close to that of Theorem
\ref{details}.
\begin{theorem}
The recursion given by the vanishing of \eqref{recurse:P2} has the unique
solution
$$
\<\tau_{k_1}\dots\tau_{k_n}\mid\lambda_g\lambda_{g-1}\>_g =
\frac{(2g+n-3)! (2g-1)!!}{(2g-1)! (2k_1-1)!!  \dots (2k_n-1)!!}
\<\tau_{g-1}\mid\lambda_g\lambda_{g-1}\>
$$
\end{theorem}

\end{document}